\documentclass[12pt]{amsart}
\usepackage{amssymb,array,epsfig}

\newtheorem{lemma}{Lemma}

\newtheorem{theorem}[lemma]{Theorem}

\newcommand{\CS}{\mathcal{S}}
\newcommand{\CH}{\mathcal{H}}

\newcommand{\CD}{\mathcal{D}}

\newcommand{\R}{\mathbb{R}}
\newcommand{\C}{\mathbb{C}}

\newcommand{\N}{\mathbb{N}}

\newcommand{\spec}{\mathrm{Spec}\,}

\newcommand{\spa}{\mathrm{Span}\,}

\newcommand{\ud}{\,\mathrm{d}}

\title[]{Finite lifetime eigenfunctions of coupled systems of
harmonic oscillators} \author[L.~Boulton,
S.A.M.~Marcantognini,
M.D.~Mor\'an]{L.~Boulton$^1$, S.A.M.~Marcantognini$^2$ and
M.D.~Mor\'an$^3$}

\date{March 2004}

\subjclass[2000]{Primary: 34L40 ; Secondary: 34L10, 81Q10.}

\keywords{Non-commutative harmonic oscillators,
higher dimensional Hermite basis, eigenfunction expansion.}

\begin{document}

\begin{abstract}
We consider a Hermite-type basis for which the eigenvalue problem
associated to the operator $H_{A,B}:=B(-\partial_x^2)+Ax^2$
acting on $L^2(\R;\C^2)$ becomes a three-terms recurrence.
Here $A$ and $B$ are $2\times 2$ constant positive definite
matrices.
Our main result provides an explicit characterization
of the eigenvectors of $H_{A,B}$
that lie in the span of the first four elements of this basis
when $AB\not= BA$.
\end{abstract}

\maketitle
\section{Introduction}
It is well known that the spectrum of the
harmonic oscillator Hamiltonian \[H_\alpha
:=-\partial_x^2+\alpha^2 x^2, \qquad \alpha>0,\]
acting on $L^2(\R)$ consists of the
non-degenerate eigenvalues \linebreak $\{\alpha (2n+1)\}_{n=0}^\infty$ with
corresponding normalized eigenfunctions \begin{equation} \label{e2}
   \phi^\alpha_n(x)=\frac{\alpha^{1/4}h_n(\alpha^{1/2} x)
   e^{-\alpha x^2/2}}{\sqrt{2^nn!\sqrt{\pi}}},
\end{equation}
where $h_n(x):=(-1)^ne^{x^2}\partial_x^n[e^{-x^2}]$ is
the $n$-th Hermite polynomial. The present paper is
devoted to studying the spectrum of a matrix version of
$H_\alpha$, the operator
\[
   H_{A,B}:=B(-\partial_x^2)+Ax^2,
\]
acting on $L^2(\R;\C^2)$, where $A$ and $B$ are two $2\times 2$
constant positive definite matrices.

In contrast to the scalar
situation,
the spectral analysis of $H_{A,B}$ is far more involved due to
the non-commutativity of the coefficients.
If $AB=BA$, it is not difficult
to find the eigenvalues and eigenfunctions of $H_{A,B}$ from those
of $H_\alpha$. On the other hand, when
$AB\not=BA$, the eigenvalues and eigenfunctions of $H_{A,B}$ are
connected to those of $H_\alpha$ in a highly non-trivial manner
(see Theorem \ref{t3} below).

Our recent interest in describing spectral properties
of operators such as $H_{A,B}$ arises from two sources.
In a series of recent works, Parmeggiani and Wakayama,
cf.
\cite{pw1}, \cite{pw2} and \cite{pw3}, characterize
the spectrum of the operator
$K_{\tilde{A},\tilde{B}}:=\tilde{A}(-
\partial_x^2+x^2)+\tilde{B}(2x\partial_x+x^2)$ acting on
$L^2(\R;\C^2)$, assuming that $\tilde{A}$ is
definite positive and $\tilde{B}=-\tilde{B}^t$.
Although the two operators are related,
it does not seem possible to obtain the eigenvalues of
$H_{A,B}$ from those of $K_{\tilde{A},\tilde{B}}$.
In \cite{pw2} and \cite{pw3}
the eigenfunctions of $K_{\tilde{A},\tilde{B}}$ are found
in terms of a twisted Hermite-type
basis of $L^2(\R;\C^2)$. In this basis the eigenvalue problem
associated to $K_{\tilde{A},\tilde{B}}$ becomes a
three-term recurrence. The strategy presented below
for analyzing the spectrum of $H_{A,B}$ will be similar.

Our second motivation is heuristic.
It is known that the scalar harmonic oscillator, $H_{\alpha}$,
achieves the optimal value for the constant
in the Lieb-Thirring inequalities with power $\sigma\geq 3/2$,
cf. \cite{lawe1}. It would be of great interest finding
Hamiltonians with similar properties for
Lieb-Thirring-type inequalities for magnetic Schr\"odinger
and Pauli
operators, cf. \cite{erso}, \cite{lawe1} and \cite{lawe2}.
Due to their close connection with the harmonic oscillator,
both $K_{\tilde{A},\tilde{B}}$ and the presently
discussed $H_{A,B}$ are strong candidates for further
investigations in this direction.

The plan of the paper is as follows. Section~2  is devoted
to describing elementary facts about $H_{A,B}$. In section~3 we
consider a
basis for which the eigenvalue
problem associated to $H_{A,B}$ is expressed as a three-term
recurrence. The main results are to be found in section~4 where
we establish necessary and sufficient conditions, given
explicitly in terms of the entries of $A$ and $B$, for an
eigenfunction of
$H_{A,B}$ to be the linear combination of the first four
elements of this basis.


\section{Elementary properties of $H_{A,B}$}

We define $H_{A,B}$ rigorously as the self-adjoint operator whose
domain, denoted below by $\CD$, is the set of all
\[
    \begin{pmatrix}
      \phi \\
      \psi
    \end{pmatrix} \in L^2(\R;\C^2)
\]
such that
\begin{equation*}
   \phi,\psi \in H^2(\R)\cap \left\{f\in L^2(\R):\int |x^2 f(x)|^2 <\infty
   \right\}=H^2(\R)\cap \widehat{H^2}(\R),
\end{equation*}
where $H^2(\R)$ denotes the Sobolev space of index (2,2)
and ``$\widehat{\ \cdot \ }$'' denotes Fourier transform.
Since
\begin{equation}\label{e1}
  \left<\!\! H_{A,B}\begin{pmatrix} \phi \\ \psi \end{pmatrix}
  \!\!,\!\!
  \begin{pmatrix} \phi \\ \psi \end{pmatrix} \!\!\right> =
  \int\! \left[B\!\begin{pmatrix} \phi' \\ \psi' \end{pmatrix}
  \!\!\cdot\!\!
  \begin{pmatrix} \phi' \\ \psi' \end{pmatrix} +
  A\!\! \begin{pmatrix} x\phi \\ x\psi \end{pmatrix}\!\!\cdot\!\!
  \begin{pmatrix} x\phi \\ x\psi \end{pmatrix}\right] \!\!\ud x \geq 0,
\end{equation}
$H_{A,B}$ is a symmetric operator.
It is well known
(cf. \cite{resiv2}) that if the domain of $H_\alpha$
is chosen to be $H^2(\R)\cap \widehat{H^2}(\R)$,
then $H_\alpha$ is self-adjoint, non-negative
and $\CS(\R)$, the Schwartz space, is a core for $H_\alpha$. Thus
$H_{A,B}$ with domain $\CD$ is a self-adjoint
non-negative operator with core $\CS(\R;\C^2)$. Indeed, these
properties are obvious when $A$ is a diagonal non-negative matrix
and $B=\mathrm{Id}$, the identity matrix. The general case follows by considering
the factorization
\begin{equation*}
   H_{A,B}=B^{1/2}E^{\ast}\,H_{C,\mathrm{Id}}\,EB^{1/2},
\end{equation*}
where $B^{-1/2}AB^{-1/2}=E^\ast C E$ is the Jordan diagonalization
of the former matrix, and by using the fact
that $\CD$ is invariant under the action of constant matrices.

\begin{lemma} \label{t1}
  The spectrum of $H_{A,B}$ consists exclusively of isolated
  eigenvalues of finite multiplicity whose only accumulation point
  is $+\infty$. Moreover, if $\lambda_n$ denotes
  the $n$-th eigenvalue of this operator counting multiplicity,
  then \[a_1^{1/2}b_1^{1/2}(2n+1)\leq \lambda_{2n+1}\leq
  \lambda_{2n+2}
  \leq a_2^{1/2}b_2^{1/2}(2n+1),\] where $0<a_1\leq a_2$ and $0<b_1\leq
  b_2$, are
  the eigenvalues of $A$ and $B$, respectively.
\end{lemma}
\begin{proof} It reduces to showing that
\[\lambda_n(H_{A_1,B_1})\leq \lambda_n(H_{A,B})
  \leq \lambda_n(H_{A_2,B_2}),\]
  where $A_j=a_j(\mathrm{Id})$ and $B_j=b_j
  (\mathrm{Id})$. This follows directly from the min-max principle
  (cf. \cite{resi}), the estimates
\begin{equation*}
    0 < \begin{pmatrix} a_1 & 0 \\0  & a_1 \end{pmatrix} \leq A
    \leq \begin{pmatrix} a_2 & 0 \\0  & a_2 \end{pmatrix} \quad
    \mathrm{and} \quad
  0 < \begin{pmatrix} b_1 & 0 \\ 0 &b_1 \end{pmatrix} \leq B
    \leq \begin{pmatrix} b_2 & 0 \\0  & b_2 \end{pmatrix},
\end{equation*}
and \eqref{e1}.
\end{proof}

The above universal bound is not sharp in general and
for most pairs $(A,B)$, $\lambda_{2n+1}\not = \lambda_{2n+2}$.

\medskip

As we mentioned earlier, it is not difficult to compute the
eigenvalues and eigenfunctions of $H_{A,B}$ when $A$ and $B$
commute. Indeed $AB=BA$ if, and only if, $A$ and $B$ have one (and
hence both) eigenvectors in common. Let $w_j\not =0$ be such that
$Aw_j=a_jw_j$ and $Bw_j=b_jw_j$, for $j=1,2$. Let
$\phi^\alpha_n(x)$ be, as in \eqref{e2}, the eigenfunctions
 of $H_\alpha$. Let $\beta_j=\sqrt{a_j/b_j}>0$. Then
\begin{align*}
  H_{A,B}w_j \phi^{\beta_j}_n(x) &= (-b_j\partial_x^2 +a_j
  x^2) w_j \phi^{\beta_j}_n(x) \\
  & = b_j(-\partial_x^2+(a_j/b_j)x^2)w_j \phi^{\beta_j}_n(x) \\
  & = b_j^{1/2}a_j^{1/2} (2n+1)w_j \phi^{\beta_j}_n(x).
\end{align*}
By choosing $\|w_j\|=1$, the family $\{w_j \phi^{\beta_j}_n(x):
j=1,2;\,n=0,1,\ldots\}$ is an orthonormal basis of $L^2(\R;\C^2)$,
hence
\[
    \spec H_{A,B}=\{b_j^{1/2}a_j^{1/2} (2n+1):j=1,2;\,n=0,1,\ldots
    \}.
\]

The analysis below will show that finding the eigenvalues and
eigenfunctions of $H_{A,B}$ whenever $AB\not=BA$ is by no means of the trivial
nature as the above case.


\section{Hermite expansion of the eigenfunctions in the
non-commutative case}

Without further mention, we will often suppress the sub-indices in
operator expressions. The structure of $H_{A,B}\equiv H$ allows us
to decompose $L^2(\R;\C^2)$ into two invariant subspaces where the
eigenvalue problem can be studied independently. We  perform
this decomposition as follows. Given $\alpha>0$, let
\begin{gather*}
\CH_+^\alpha := \spa \{v x^je^{-\alpha x^2/2}: v\in \C^2,j=2k,
k=0,1,\cdots \},
   \\
\CH_-^\alpha := \spa \{v x^je^{-\alpha x^2/2}: v\in\C^2,j=2k+1,
k=0,1,\cdots\},
\end{gather*}
and denote by $H^\pm=H|(\CD\cap\CH^{\alpha}_{\pm})$.
Since $\CH_\pm^\alpha$ are invariant under $\partial_x^2$,
multiplication by $x^2$ and action of constant
matrices, these
spaces are also invariant under $H$. Hence
$H^\pm:\CD \cap
\CH_\pm^\alpha\longrightarrow \CH_\pm^\alpha$ are self-adjoint operators
and
\[
    \spec H=\spec H^+ \cup
    \spec H^-.
\]

Let
\[
   L=2^{-1/2}(x+\partial_x) \quad \mathrm{and} \quad
   L^\ast =2^{-1/2} (x-\partial_x)
\]
be the annihilation and creation operators for the scalar harmonic
oscillator. Then
\[
   L\phi^1_0=0,\quad L\phi^1_n=n^{1/2}\phi^1_{n-1} \quad
 \mathrm{and} \quad L^\ast \phi^1_n=(n+1)^{1/2}\phi^1_{n+1}.
\]
From these relations one can easily deduce the recurrent identities
\begin{gather*}
  2\alpha x^2 \phi^\alpha_n\!=\!(n+2)^{1/2}(n+1)^{1/2}
  \phi^\alpha_{n+2}
  +(2n+1)\phi_n^\alpha+n^{1/2}(n-1)^{1/2}\phi^\alpha_{n-2}, \\
  2\alpha^{-1}\partial_x^2 \phi^\alpha_n\!=\!(n+2)^{1/2}(n+1)^{1/2}
  \phi^\alpha _{n+2}
  \!-\!(2n+1)\phi_n^\alpha\!+\!n^{1/2}(n-1)^{1/2}\phi^\alpha_{n-2},
\end{gather*}
where, here and elsewhere, any quantity with negative sub-index is
zero. Since $\{\phi^\alpha_n\}_{n=0}^\infty$ is an orthonormal
basis for $L^2(\R)$, we can expand any vector of $L^2(\R;\C^2)$
via
\begin{equation} \label{e3}
   \begin{pmatrix}\phi \\ \psi \end{pmatrix}=
   \sum_{n=0}^\infty v_n \phi_n^\alpha,
\end{equation}
for a suitable unique sequence $(v_n)\in l^2(\N;\C^2)$. Moreover,
denoting by $N_\alpha:=\alpha^{-1}A-\alpha B$ and
 $M_\alpha:=\alpha^{-1}A+\alpha B$,
\begin{align*}
   H\begin{pmatrix}\phi \\ \psi \end{pmatrix}&=\left[(-\partial_x^2)B+x^2A
\right]\sum_{n=0}^\infty v_n \phi_n^\alpha \\
  &=\sum_{n=0}B v_n (-\partial_x^2)
  \phi_n^\alpha+Av_nx^2\phi_n^\alpha \\
  &=\frac{1}{2} \sum_{n=0}(n+2)^{1/2}(n+1)^{1/2}N_\alpha v_n\phi_{n+2}^\alpha
  +(2n+1)M_\alpha v_n \phi_n^\alpha + \\ & \quad
  +n^{1/2}(n-1)^{1/2} N_\alpha v_n
 \phi_{n-2}^\alpha \\
 &=\frac{1}{2} \sum_{k=0}\left[ k^{1/2}(k-1)^{1/2}N_\alpha v_{k-2}
 +(2k+1)M_\alpha v_k+ \right.\\ &  \quad \left.+(k+2)^{1/2}(k+1)^{1/2}
  N_\alpha v_{k+2}\right]\phi_k^\alpha .
\end{align*}
Thus $2H^\pm$ are, respectively, similar to the block tri-diagonal
matrices
\begin{equation} \label{e6}
   \begin{pmatrix}
  S^\pm_0 & T^\pm_1 &     &            \\
  T^\pm_1 & S^\pm_1 & T^\pm_2 &            \\
      & T^\pm_2 & S^\pm_2 & T^\pm_3        \\
      &     & T^\pm_3  & \ddots
  \end{pmatrix}
\end{equation}
acting on $l^2(\N;\C^2)$, where
\begin{gather*}
   S^+_k=(4k+1)M_\alpha, \qquad T_k^+=(2k)^{1/2}(2k-1)^{1/2}N_\alpha, \\
   S^-_k=(4k+3)M_\alpha \quad \mathrm{and} \quad
   T_k^-=(2k)^{1/2}(2k+1)^{1/2}N_\alpha .
\end{gather*}

\medskip

In order to reduce the amount of notation in our subsequent
discussion, we consider $H_{A,B}$ in canonical form as follows.
If $0\leq b_1\leq b_2$ are the eigenvalues of
$B$, let $U^\ast( \mathrm{diag}[b_1,b_2]) U$ be the diagonalization
of $B$, and set $\tilde{A}:= b_1^{-1}UAU^\ast$ and
$\tilde{B}:=\mathrm{diag}[1,b_2/b_1]$. Then
\[
   H_{A,B}=b_1U^\ast
   H_{\tilde{A},\tilde{B}}U
\]
where
\begin{equation} \label{e12}
   \tilde{A}=\begin{pmatrix} a & \xi \\ \overline{\xi} & c \end{pmatrix}
   \qquad \mathrm{and} \qquad \tilde{B}=\begin{pmatrix} 1 & 0 \\ 0 & b
   \end{pmatrix}.
\end{equation}
Here the positivity of $A$ and $B$ is equivalent to the conditions
\[
    b\geq 1,\qquad a,c>0\qquad \mathrm{and}\qquad
    0\leq|\xi|^2<ac.
\]
Furthermore notice that $AB=BA$ if, and only if, either $b=1$ or
$\xi=0$. Hence, unless otherwise specified, we will consider
without loss of generality that the pair $(A,B)$ is always the pair
$(\tilde{A},\tilde{B})$ in \eqref{e12}.

\medskip

By virtue of the tri-diagonal representation \eqref{e6},
it seems natural to expect that the Hermite series \eqref{e3} may be
a good candidate
for expanding the eigenfunctions of $H$. Not to mention that it is
the obvious extension of the scalar and commutative cases.
In this respect, we may consider ``finite lifetime'' series
expansions of eigenfunctions $\Phi$ of $H$,
\begin{equation} \label{e4}
     \Phi=\sum _{n=0}^{k} v_n \phi_n^{\alpha(k)}
\end{equation}
for suitable finite $k\in \N\cup\{0\}$, $\alpha(k)>0$ and $v_n\in
\C^2$. The results we present below show that, contrary to the
above presumption, \eqref{e3} is not such a
good candidate for expanding $\Phi$
for small values of $k$. To be more precise, we show that
for $k=0,1,2,3$, an expansion of type
\eqref{e4} is allowed only for a small sub-manifold of the region
\begin{equation} \label{e8}
   R:=\{(b,a,c,|\xi|)\in \R^4\,:\, a,c>0,\, b\geq 1,\, 0\leq |\xi|^2<ac\}
\end{equation}
corresponding to all positive definite pairs $(A,B)$.

We first discuss the cases $k=0,1$ and leave $k=2,3$
for the forthcoming section. The following result
includes a family of test bases larger than the one
considered in \eqref{e4}.

\begin{lemma} \label{t2}
Let $\phi^\alpha_n(x)$ be the eigenfunctions of the scalar
harmonic oscillator $H_\alpha$. Then $H_{A,B}$ has an
eigenfunction of the type
$
   \Phi(x)=(\tilde{a}\phi^{\alpha}_n(x) , \tilde{b}\phi^\beta_m(x)
   )^\mathrm{t}
$
where $\tilde{a},\tilde{b}\in \C$, $\alpha,\beta>0$ and
$m,n\in \N\cup \{0\}$, if and only if $AB=BA$.
\end{lemma}

\begin{proof}
If $H\Phi=\lambda \Phi$, then
\begin{gather*}
-\tilde{a}(\phi_n^\alpha)''+a x^2\tilde{a}\phi^\alpha_n+\xi
x^2\tilde{b}\phi^\beta_m-\lambda \tilde{a}\phi^\alpha_n=0, \\
-b\tilde{b}(\phi_m^\beta)''+c
x^2\tilde{b}\phi^\beta_m+\overline{\xi}
x^2\tilde{a}\phi^\alpha_n-\lambda \tilde{b}\phi^\beta_m=0.
\end{gather*}
If $\tilde{a}=0$ or $\tilde{b}=0$ in the above identities,
necessarily
$\xi=0$ so $AB=BA$. Hence without loss of generality
we can assume that
$\tilde{a}\tilde{b}\not=0$.

If $\alpha\not=\beta$, once again
$\xi=0$. Then we may suppose that $\alpha=\beta$.
Since both left hand sides of the above identities are equal to
$p(x)e^{-\alpha x^2/2}$, where in both cases $p(x)$ is a
polynomial of degree $2+\mathrm{max}(m,n)$, necessarily either
$\xi=0$ or $m=n$. In the latter case, the above system is
rewritten as
\begin{gather*}
-\tilde{a}(\phi_n^\alpha)''+(a\tilde{a}+\xi\tilde{b})
x^2\phi^\alpha_n-\lambda \tilde{a}\phi^\alpha_n=0, \\
-b\tilde{b}(\phi_n^\alpha)''+(c\tilde{b}+\overline{\xi}\tilde{a})
x^2\phi^\alpha_n-\lambda \tilde{b}\phi^\alpha_n=0.
\end{gather*}
Since $\phi_n^\alpha$ is an eigenfunction of $H_\alpha$ where
$\alpha>0$,
necessarily $\xi\in \R$. Furthermore,
\begin{gather*}
a+\xi \tilde{b}/\tilde{a}=\alpha^2=c/b+\xi\tilde{a}/(b\tilde{b})
\qquad \mathrm{and}\\
a+\xi\tilde{b}/\tilde{a}=\lambda/(2n+1)=bc+\xi\tilde{a}b/\tilde{b}.
\end{gather*}
Hence necessarily $b=1$.
\end{proof}

Since $\CH^\pm$ are invariant under the action of $H$, and the
even (resp. odd) terms in the series \eqref{e4} belong to $\CH^+$
(resp. $\CH^-$), the above lemma ensures that
$
  \Phi=v_0\phi_0^\alpha+v_1\phi_1^\alpha
$
is an eigenfunction of $H$ if, and only if, $A$ and $B$ commute.


\section{Four-term expansion of eigenfunctions of $H_{A,B}$}
In this section we study necessary and sufficient conditions in
order to guarantee that $\Phi\in L^2(\R;\C^2)$, with finite
lifetime expansion of the type \eqref{e4} for $k=2$ and 3, is an
eigenfunction of $H$ for suitable $\alpha(k)>0$ when $AB\not=BA$.
In other words, assuming that $\Phi$ satisfies the constraint
\[
   \Phi=v_0\phi_0^\alpha+v_1\phi_1^\alpha+v_2\phi_2^\alpha+
   v_3\phi_3^\alpha,
\]
we aim to investigate conditions ensuring
$H\Phi=\lambda\Phi$.

Since
\[v\phi_{2n}^\alpha\in \CH_+^\alpha \quad \mathrm{and} \quad
v\phi_{2n+1}^\alpha\in \CH_-^\alpha, \quad v\in \C^2\] for all
$n\in \N\cup\{0\}$, and the subspaces $\CH_\pm^\alpha$ are invariant
under $H$, we may consider the even and odd cases separately. To
this end, let
\begin{equation} \label{e5}
 \Phi^+=v_0\phi_0^\alpha+v_2\phi_2^\alpha \quad \mathrm{and}
 \quad \Phi^-=v_1\phi_1^\alpha+v_3\phi_3^\alpha,
\end{equation}
$\Phi^\pm\in H^\pm$ respectively. Then our goal is to find
necessary and sufficient conditions, given in terms of
$(b,a,c,|\xi|)\in R\setminus \partial R$, ensuring that $\Phi^\pm$ is an
eigenfunction of $H^\pm$. The following is our main result.

\begin{theorem} \label{t3}
Let
\[
    B=\begin{pmatrix} 1 &0 \\ 0 & b \end{pmatrix} \qquad \mathrm{and}
    \qquad A=\begin{pmatrix} a &\xi \\ \overline{\xi} & c
    \end{pmatrix},
\]
where $b>1$, $a,c>0$ and $0<|\xi|^2<ac$. Let $\beta>0$ be such
that $\beta^2$ is an eigenvalue of $B^{-1/2}AB^{-1/2}>0$. Let
\begin{equation*}
\lambda_{\mathrm{even}}:=
\frac{5\beta(ab+c-2\beta^2b)}{a+c-(b+1)\beta^2}
\quad \mathrm{and} \quad
\lambda_\mathrm{odd}:= \frac{7\beta(ab+c-2\beta^2b)}{a+c-(b+1)
    \beta^2}.
\end{equation*}
Then
\item[i)] $H^+$ has an eigenfunction $\Phi^+(x)$ of type \eqref{e5} if
and only if
\begin{equation} \label{e11}
   2\lambda_\mathrm{even}[a+c-(b+1)\beta^2]=5\beta(\beta-
   \lambda_\mathrm{even})
   (\lambda_\mathrm{even}-\beta b).
\end{equation}
In this case $\alpha=\beta$ and $H^+\Phi^+=\lambda_\mathrm{even}
\Phi^+$. Furthermore, $\lambda_\mathrm{even}$ is an eigenvalue
of $M_\beta$.
\item[ii)]  $H^-$ has an eigenfunction
$\Phi^-(x)$ of type \eqref{e5} if
and only if
\begin{equation} \label{e13}
   6\lambda_\mathrm{odd}[a+c-(b+1)\beta^2]=7\beta(
   3\beta-\lambda_\mathrm{odd})
   (\lambda_\mathrm{odd}-3\beta b).
\end{equation}
In this case $\alpha=\beta$ and $H^-\Phi^-=\lambda_\mathrm{odd}
\Phi^-$. Furthermore $\lambda_\mathrm{odd}$ is an eigenvalue of
$3M_\beta$.
\end{theorem}

Notice that the conditions on $a,b,c$ and $\xi$ ensure that
$AB\not =BA$.

\begin{proof}
Put
\[
   \Phi^+(x)=(u_0+u_2x^2)e^{-\alpha x^2/2},
\]
for $u_0,u_2\not=0$. Since $H^+$ is similar to the tri-diagonal
matrix \eqref{e6}, then $H^+\Phi^+=\lambda \Phi^+$, if and only if
\begin{align*}
   & S_0^+u_0+T^+_1u_2=\lambda u_0 \\
   & T^+_1u_0+S_1^+u_2=\lambda u_2 \\
   & T^+_2 u_2=0.
\end{align*}
The latter equation implies that $u_2\in \ker N_\alpha$ and thus
the first one implies that $\lambda$ is an eigenvalue of
$S_0^+=M_\alpha$ with associated eigenfunction $u_0$. A
straightforward computation shows that the above system is
equivalent to
\begin{equation} \label{e7}
\begin{aligned}
   & (A-\alpha^2B)u_2=0 \\
   & (A-\alpha^2B)u_0+(5\alpha B-\lambda)u_2=0 \\
   & (\alpha B-\lambda)u_0-2Bu_2=0.
\end{aligned}
\end{equation}
The first equation holds if and only if $\alpha=\beta$. Here
\begin{equation} \label{e14}
   \beta=+\sqrt{\frac{ab+c\pm\sqrt{(c-ab)^2+4|\xi|^2b}}{2b}}
\end{equation}
and
$
    u_2=\begin{pmatrix} c-\beta^2b \\ -\overline{\xi}
    \end{pmatrix}.
$
Notice that in this case
\[
  0=\det (A-\beta^2B)= (a-\beta^2)(c-\beta^2 b)-|\xi|^2.
\]
 Let
$
    \tilde{u}_2=\begin{pmatrix} a-\beta^2 \\ \xi
    \end{pmatrix}.
$
Then $\tilde{u}_2\perp u_2$ and
\[
   (A-\beta^2B)\tilde{u}_2=(a+c-(1+b)\beta^2) \tilde{u}_2.
\]
Decompose
\[
   u_0=\gamma u_2+\tilde{\gamma} \tilde{u}_2,
\]
for suitable $\gamma,\tilde{\gamma}\in \C$. Then the second
identity of \eqref{e7} holds if and only if,
$\lambda=\lambda_\mathrm{even}$  and
\begin{equation} \label{e10}
   \tilde{\gamma}(a+c-(1+b)\beta^2)^2(a-\beta^2)=5\beta(1-b)|\xi|^2.
\end{equation}
The third identity of \eqref{e7} can be rewritten as the system
\begin{align*}
   &
   \gamma(\beta-\lambda)(c-\alpha^2b)+\tilde{\gamma}(\beta-\lambda)
   (a-\beta^2)=2(c-\beta^2 b) \\
   &\gamma(\beta b-\lambda)(-\overline{\xi})+\tilde{\gamma}
   (\beta b-\lambda) \overline(\xi)=-2b\overline{\xi}
\end{align*}
in $\gamma$ and $\tilde{\gamma}$. By finding $\tilde{\gamma}$ from
this system (for instance by Newton's method) and by equating to
\eqref{e10}, a straightforward computation yields \eqref{e11}.

The proof of ii) is similar.
\end{proof}

\medskip

Notice that there is a duality of signs in the definition of
$\beta$ (cf. \eqref{e14}), so conditions \eqref{e11} and
\eqref{e13} comprise two possibilities each. Let
\[
    \Omega_\mathrm{even}^\pm:=\{(b,a,c,|\xi|)\in R\,:\,
    (\ref{e11}) \mathrm{\ holds}\}
\]
and
\[
    \Omega_\mathrm{odd}^\pm:=\{(b,a,c,|\xi|)\in R\,:\,
    (\ref{e13}) \mathrm{\ holds}\},
\]
where the sign for the super-index is chosen in concordance to the
sign in expression \eqref{e14}. By computing the partial derivatives of
both sides of identities \eqref{e11} and \eqref{e13}, a
straightforward but rather long computation which we omit in the
present discussion, shows that these four regions are smooth
3-manifolds embedded in $R\subset \R^4$, see \eqref{e8}.

The fact that
$\Omega_\mathrm{even}^\pm$ are non empty is consequence of the
following observation. By fixing $\tilde{b}>1$ and
$|\tilde{\xi}|>0$, and by putting $c=a\tilde{b}$, condition
\eqref{e11} can be rewritten as
\begin{align*}
0&=\tilde{b}^{-1}(\tilde{b}+1)^2[2(a\tilde{b}+c-2\beta
^2\tilde{b})-(\beta-\lambda)(\lambda-\tilde{b}\beta)]
\\
  &= \pm \frac{4|\tilde{\xi}|}{\sqrt{\tilde{b}}}(\tilde{b}+1)^2-\beta^2
  (1-9\tilde{b})(9-\tilde{b}),
\end{align*}
where $\beta^2=a\pm|\tilde{\xi}|\tilde{b}^{-1/2}$. Then
$(\tilde{b},\tilde{a},\tilde{a}\tilde{b},|\tilde{\xi}|)\in
\Omega^+_\mathrm{even}$ whenever
\[
   \tilde{a}=-\frac{|\tilde{\xi}|(5\tilde{b}^2-90\tilde{b}+5)}
   {\tilde{b}^{1/2}(9\tilde{b}^2-82\tilde{b}+9)}
\]
for $9<\tilde{b}<9+4\sqrt{5}$, and
$(\tilde{b},\tilde{a},\tilde{a}\tilde{b},|\tilde{\xi}|)\in
\Omega^-_\mathrm{even}$ whenever
\[
   \tilde{a}=\frac{|\tilde{\xi}|(5\tilde{b}^2-90\tilde{b}+5)}
   {\tilde{b}^{1/2}(9\tilde{b}^2-82\tilde{b}+9)}
\]
for $1<\tilde{b}<9$ or $\tilde{b}>9+4\sqrt{5}$. Furthermore,
notice that if
$(\tilde{b},\tilde{a},\tilde{a}\tilde{b},|\tilde{\xi}|)\in
\Omega_\mathrm{even}^\pm$, then
\[
   (\tilde{b},r\tilde{a},r\tilde{c},r|\tilde{\xi}|)
   \in \Omega_\mathrm{even}^\pm, \qquad \mathrm{for\ all}\quad
   r>0.
\]
This shows that $\Omega_\mathrm{even}^\pm$ are non-empty,
 unbounded and the semi-axis $\{(b,0,0,0):b>1\}$ intersects
$\partial\Omega_\mathrm{even}^\pm$. All these properties
also hold for $\Omega_\mathrm{odd}^\pm$.

Furthermore,
\[
   \Omega^+_\mathrm{even} \cap
   \Omega^+_\mathrm{odd}=\varnothing\qquad \mathrm{and}\qquad
   \Omega^-_\mathrm{even} \cap
   \Omega^-_\mathrm{odd}=\varnothing.
\]
Indeed, if
\eqref{e11} and \eqref{e13} hold at the
same time for the same tetrad $(b,a,c,|\xi|)$, then, according to
Theorem~\ref{t3},  $\lambda$ should be at the same
time eigenvalue of $M_\beta$ and $3M_\beta$ ($\beta$ defined with
the same sign). Obviously the latter is a contradiction.
\begin{figure}[t]
\epsfig{file=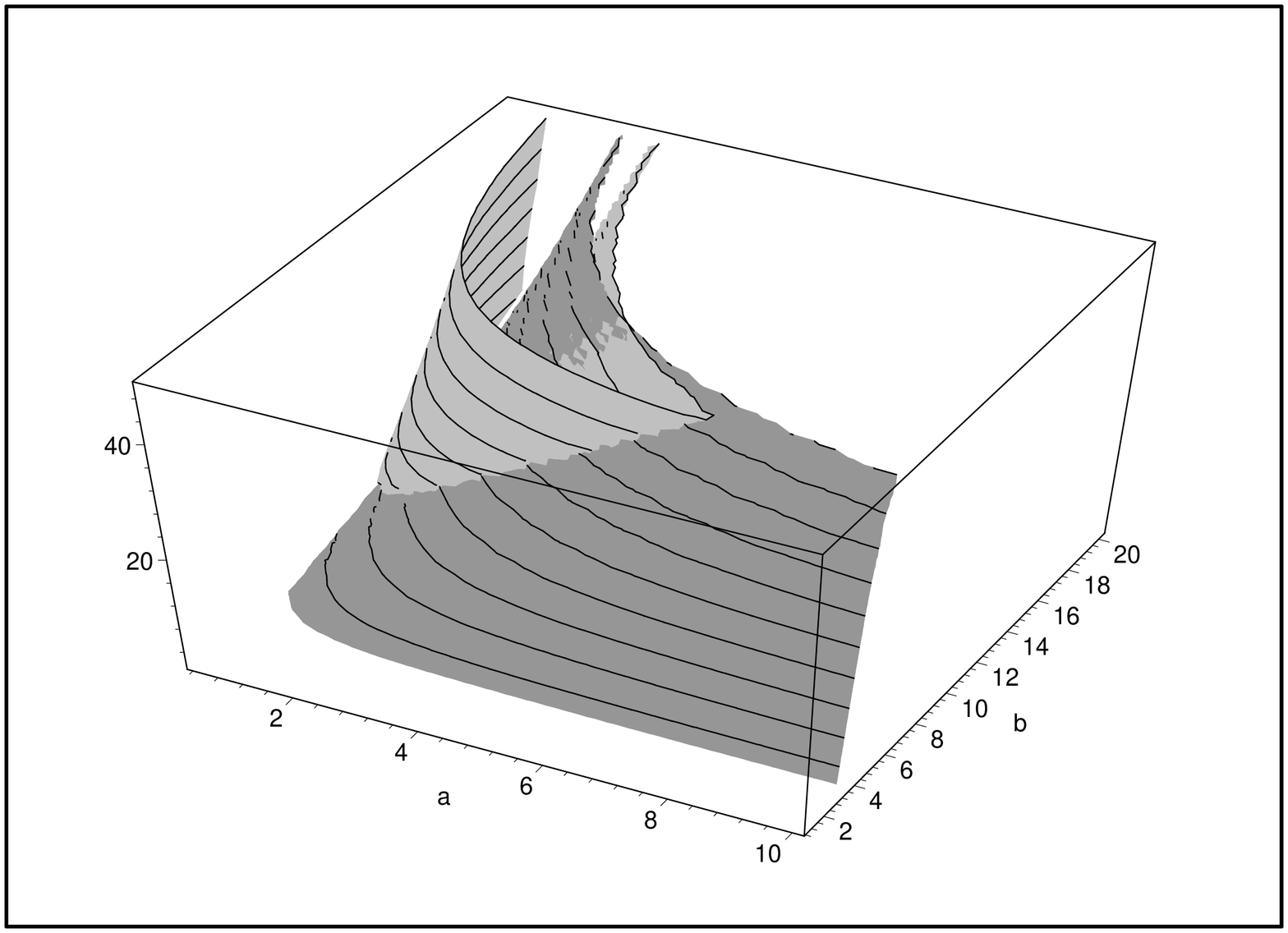, angle=-90, width=2.3in} \quad
\epsfig{file=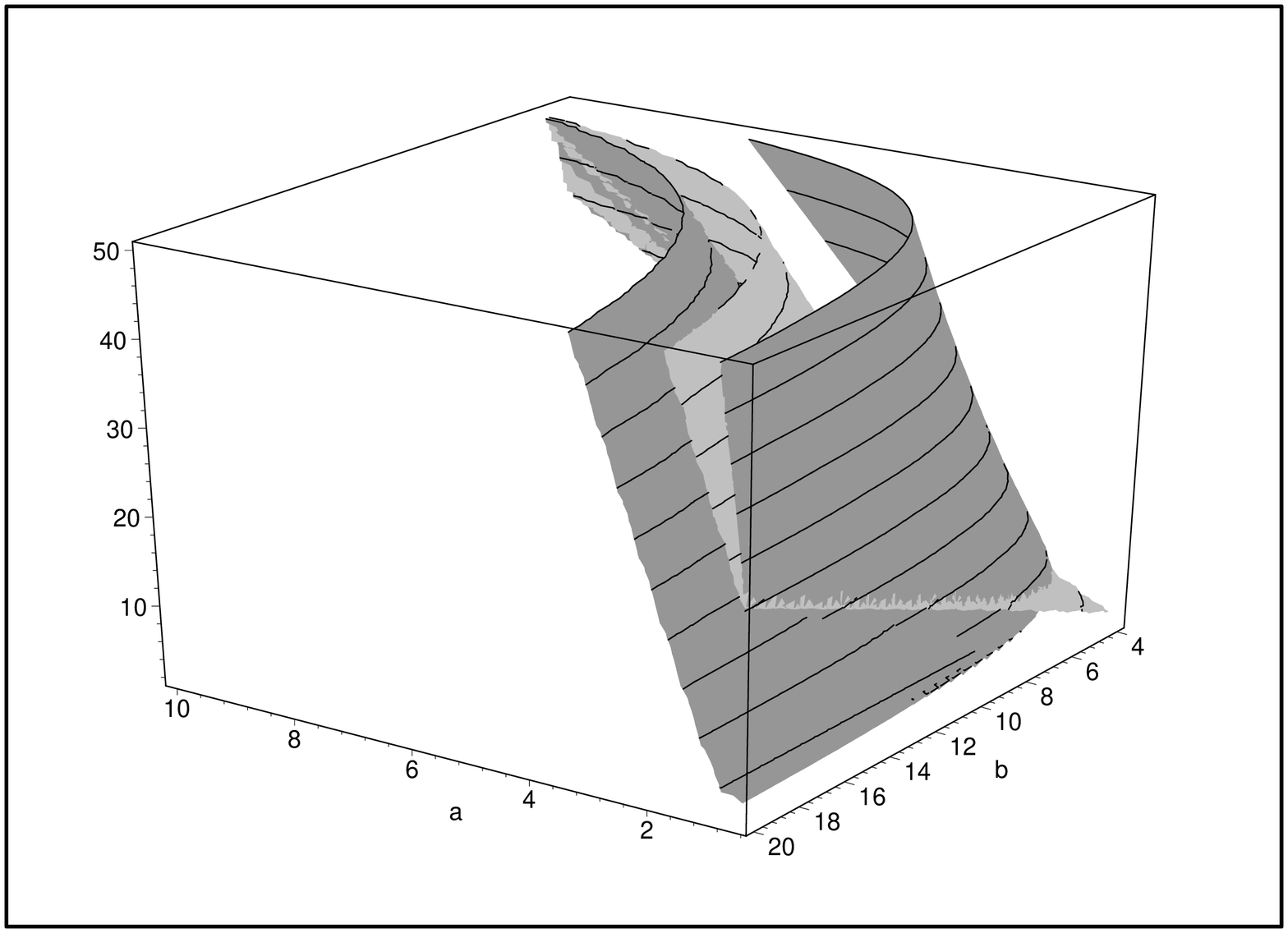, angle=-90, width=2.3in}
\caption{Left: projection
of $\Omega_\mathrm{even}^+$ and $\Omega_\mathrm{even}^-$ onto the
hyperplane $|\xi|=1$, the darker color corresponds to
$\Omega_\mathrm{even}^+$. Right: projection
of $\Omega_\mathrm{odd}^+$ and $\Omega_\mathrm{odd}^-$ onto the
hyperplane $|\xi|=1$. The lighter color corresponds to
$\Omega_\mathrm{odd}^+$.}
\end{figure}
In figure~1 we reproduce the
projections of these four regions, onto the hyper-plane $|\xi|=1$.
These picture suggest that $\Omega^+\cap \Omega^-\not =
\varnothing$.

\medskip

Finally, we may comment on the issue of considering a more general
basis for expanding the eigenfunctions of $H$. One might think
that a natural candidate
for generalizing \eqref{e4} is the finite expansion
\[
   \Phi=\sum_{j=0}^m\begin{pmatrix} a_{j}\phi_j^{\alpha} \\
    b_{j}\phi_j^{\beta}
    \end{pmatrix},
\]
where $a_{j},b_{j}$ are complex number and $\alpha,\beta>0$. We
studied a particular case of this in Lemma~\ref{t2}. It turns out
that if $\Phi$ is an eigenfunction of the above form, then either
$\xi=0$ or $\beta=\alpha$, so it should be as in \eqref{e4}. This
can be easily proven by writing down the system for the eigenvalue
equation and considering the asymptotic behaviour of the
identities as $x\to \infty$.



\vspace{1.5in}

\begin{minipage}{2.3in}
{\scshape $^1$Lyonell Boulton}\\
{\footnotesize Department of Mathematics\\ and
Statistics, \\ University of Calgary, \\
2500 University Drive,\\
Calgary, AB, Canada T2N 1N4 \\
email: \texttt{lboulton@math.ucalgary.ca}}
\end{minipage}
\quad
\begin{minipage}{2.3in}
{\scshape $^2$Stefania Marcantognini}\\
{\footnotesize Departamento de Matem\'aticas, \\
Instituto Venezolano de \\ Investigaciones Cient\'\i ficas, \\
Apartado 21827,\\
Caracas, 1020A, Venezuela.\\ email:
\texttt{smarcant@ivic.ve}}
\end{minipage}

\bigskip

\begin{minipage}{2.3in}
{\scshape $^3$Mar\'\i a Dolores Mor\'an}\\
{\footnotesize Escuela de Matem\'aticas, \\
Facultad de Ciencias, \\
Universidad Central de Venezuela, \\
Apartado 20513,\\
Caracas, 1020A, Venezuela.\\ email:
\texttt{mmoran@euler.ciens.ucv.ve}}
\end{minipage}


\begin{thebibliography}{10}
\bibitem{erso}{\scshape L.~Erd\"os, J.P.~Solovej}
``Magnetic Lieb-Thirring inequalities with optimal
dependence on the field strength''. Preprint 2003,
\textsf{arXiv:math-ph/0306066}.
\bibitem{lawe1}{\scshape A.~Laptev, T.~Weidl},
``Recent results on Lieb-Thirring inequality'',
\emph{Universit\'e de Nantes. Exp.} XX (2000) 1-14.
\bibitem{lawe2}{\scshape A.~Laptev, T.~Weidl},
``Sharp Lieb-Thirring inequalities in high dimensions'',
\emph{Acta Math.} 184 (2000) 87-111.
\bibitem{pw1}{\scshape A.~Parmeggiani, M.~Wakayama},
``Oscillator representations and systems of ordinary differential
equations'', \emph{Proc. Natl. Acad. Sci. USA} 98
(2001) 26-30.
\bibitem{pw2}{\scshape A.~Parmeggiani, M.~Wakayama},
``Non-commutative harmonic oscillators I'', \emph{Forum Math.} 14
(2002) 538-604.
\bibitem{pw3}{\scshape A.~Parmeggiani, M.~Wakayama},
``Non-commutative harmonic oscillators II'', \emph{Forum Math.} 14
(2002) 669-690.
\bibitem{resiv2}{\scshape M. Reed, B. Simon}, {\em Methods of modern
mathematical physics, volume 2: self-adjointness}, Academic press,
New York, 1975.
\bibitem{resi}{\scshape M. Reed, B. Simon}, {\em Methods of modern
mathematical physics, volume 4: analysis of operators}, Academic
press, New York, 1978.
\end{thebibliography}
\end{document}